\documentclass[12pt]{article}
\usepackage{amssymb, amsmath, amscd, amsthm, setspace, latexsym}
\usepackage[toc,page]{appendix}
\title{Some formulas for the smallest number of generators for finite direct sums of matrix algebras}
\author{R.V.~Kravchenko\\
\small{Department of Mathematics,}\\
\small{Texas A\&M University, College Station, TX 77843-3368, USA.}\\
\small{E-mail: rkchenko@math.tamu.edu}
\and
B.V.~Petrenko\\
\small{Department of Mathematics,}\\
\small{SUNY Brockport, 350 New Campus Drive}\\
\small{E-mail: bpetrenk@brockport.edu}}
\newtheorem{theorem}{Theorem}[section]
\newtheorem{lemma}[theorem]{Lemma}

\newtheorem{corollary}[theorem]{Corollary}

\newtheorem{problem}[theorem]{Problem}

\newtheorem{definition}[theorem]{Definition}

\newtheorem{facts}[theorem]{Facts}

\newtheorem{observation}[theorem]{Observation}
\def\gen{\text{gen}}
\def\tr{\text{tr}}
\begin{document}
\maketitle

\begin{abstract} We obtain an asymptotic upper bound for the smallest number
of generators for a finite direct sum of matrix algebras with entries in a
finite field.
This produces an upper bound for a similar quantity for integer matrix rings.
We also obtain an exact formula for the smallest number
of generators for a finite direct sum of $2$-by-$2$ matrix algebras with entries in a
finite field and as a consequence obtain a formula for a similar quantity for
a finite direct sum of $2$-by-$2$ integer matrix rings. We remark that a generating set
the ring $\bigoplus_{i=1}^k M_{n_i}(\mathbb{Z})^{n_i}$ may be used as a generating set of any matrix algebra $\bigoplus_{i=1}^k M_{n_i}(R)^{n_i}$ where $R$ is an associative ring with a two-sided 1.
\end{abstract}

2000 MSC: 16S50, 15A30, 15A33, 15A36, 15A30, 16P90

\tableofcontents

\section{Introduction}
{\em All rings and algebras in this paper are
assumed associative with a two-sided identity element. As usual, in a direct sum of rings or algebras the operations are defined componentwise. When we discuss the ring $M_n(R)$
of $n$-by-$n$ matrices with entries in a ring $R$, we assume that $n \geq 2$ unless noted otherwise.}
Our paper has been inspired by the work of Philip Hall \cite{Hall}, where he studied
the smallest number of generators needed for finite direct products of various finite groups.
In particular, he showed that a direct product of up to $19$ copies of the alternating group on $5$ symbols
can be generated by $2$ elements, but not the direct product of $20$ of copies of this group.
In this paper, in particular, we provide the following formula for the smallest number of generators for a finite direct sum of the ring $M_2(\mathbb{Z})$ with itself. Given an integer $m \geq 2$, let $g(m)$ denote the largest integer with the property that
a direct sum of $g(m)$ copies of the ring $M_2(\mathbb{Z})$ has $m$ generators, then
\begin{equation}\label{salvo}
g(m) = \frac{16^m - 3 \cdot 8^{m} + 2 \cdot 4^m}{6}\text{.}
\end{equation}
In particular, $g(2) = 16$ which means that a direct sum of up to $16$ copies of the matrix
ring $M_2(\mathbb{Z})$ of $2$-by-$2$ integer matrices can be generated by $2$ elements,
while the smallest number of generators for the ring $M_2(\mathbb{Z})^{17}$ is $3$.
This result can be extended to an arbitrary ring as follows.
Suppose that $a,b$ generate the ring $M_2(\mathbb{Z})^{16}$, $R$ is a ring, and
$\phi: M_2(\mathbb{Z})^{16} \to M_2(R)^{16}$ a ring homomorphism induced by the ring homomorphism $\mathbb{Z} \to R$. We see that any element of $M_2(R)^{16}$
can be written as a sum of the terms $rw$, where $r \in R$ and $w$ is a word in
$\phi(a)$ and $\phi(b)$ (an example of such a word is $\phi(b)\phi(a)^2 \phi(b)\phi(a)$).

The paper of Petrenko and Sidki \cite[Theorem 3.11 (2)]{PetSid} contains a presentation showing that any finite direct sum of matrix algebras with entries in an infinite field
always admits two generators (see Appendix \ref{Kemer} for details).
This is no longer true in general for finite direct sums of matrix rings with entries in a finite field. Therefore, the same conclusion applies to
finite direct sums of integer matrix rings.
\footnote{ Nevertheless, by Theorem 3.11 (4) of \cite{PetSid} this is true for a finite direct sum of integer matrix rings such that the rings of the same size appear no more than three times. The number $3$ in the previous sentence is not optimal.}
Consider, for example, a direct sum $M_n(\mathbb{Z})^m$ of $m$ copies of the ring $M_n(\mathbb{Z})$. If for any $m$ this ring had $2$ generators, then the same would be true of its epimorphic image
$M_n(\mathbb{F}_2)^m$, where $\mathbb{F}_2$ is a field with two elements. Let $a = (A_1, \ldots, A_m)$ and $b = (B_1, \ldots, B_m)$ be generators of $M_n \left(\mathbb{F}_2 \right)^m$. Since $M_n(\mathbb{F}_2)$ is a finite set,
it follows that if $m$ is sufficiently large, then there exist $i\neq j$ such that $A_i = A_j$ and
$B_i= B_j$. Therefore, $a$ and $b$ cannot generate $M_n \left(\mathbb{F}_2 \right)^m$ because the ring generated by $\left(A_i, A_j\right)$ and $\left(B_i, B_j\right)$ is isomorphic to a subring of $M_n \left(\mathbb{F}_2 \right)$, and not to the larger ring
$M_n \left(\mathbb{F}_2 \right)^2$ as it should.
Therefore, if $m \to \infty$, then so is the smallest number of generators of the rings $M_n \left(\mathbb{F}_2 \right)^m$ and
$M_n(\mathbb{Z})^m$.

The more general question about the smallest number of generators for a finite direct sum of matrix
algebras of different sizes reduces to the same question about a finite direct sum of copies of the same
algebra, in view of Theorem \ref{th6}. It states that the smallest number of generators of a finite direct sum of matrix algebras is the maximum of the smallest number of generators in all the sums of terms of the same size.

\begin{definition} Let $R$ be a commutative ring. We introduce the sequence
$ \{ \gen_{m,n}(R) \}$ defined by the property that $\gen_{m,n}(R)$ is the integer such that
the $R$-algebra $M_n(R)^{\gen_{m,n}(R)}$ admits $m $ generators, while $M_n(R)^{1+\gen_{m,n}(R)}$
needs at least $m+1$ generators. If any $R$-algebra $M_n(R)^k$ admits $\leq m$ generators,
then put $\gen_{m,n}(R) = \infty$.
\end{definition}
In particular, $g(m)$ in (\ref{salvo}) is $\gen_{m,2} (\mathbb{Z})$; $\gen_{2,n}(F) = 2$ for
any $n \geq 2$ and an infinite field $F$ by
\cite[Theorem 3.11 (2)]{PetSid}; $\gen_{1,1}(F) = \infty$ for any infinite field $F$
by the formula for the Vandermonde determinant.

\medskip

Theorem \ref{th3} below implies that a set $S$ generates a finite direct sum of integer matrix rings if and only if
the reduction of $S$ modulo every prime $p$ generates the corresponding finite direct sum of matrix rings over $\mathbb{F}_p$.

This result is quantified by Theorem \ref{Lenstra1} that together with proof has been kindly communicated to us
by Hendrik Lenstra \cite{Lenstra_C}. The importance of this result to us becomes obvious if one looks at the evolution of our preprint arXiv:math/0611674.
\begin{theorem}\label{Lenstra1}
Let  $R$  be a ring whose additive group is finitely generated.
For each prime $p$,  let
            \[r(p)\]
be the smallest number of generators of
$R/pR$ as a ring, and let
            \[r(0)\]
be the smallest number of generators of
$ R\otimes_{\mathbb{Z}} \mathbb{Q}$ as a $\mathbb{Q}$-algebra. Finally, let
            \[r\]
be the smallest number of generators of  $R$  as a ring. Then
\begin{enumerate}
\item For each  prime number  $p$  we have $r(0) \le r(p) \le r$.
\item For all but finitely many prime numbers we have $r(0) = r(p)$.
\item If there is a prime number  $p$  such that  $r(0) < r(p)$,  then
                        \[r = \max\{r(p): p~ \text{prime}\}.\]
\item If $r(0) = r(p)$  for all  $p$,  then
            either  $r = r(0)$  or  $r = r(0) + 1$.
\end{enumerate}
\end{theorem}

This theorem is later restated as Theorem \ref{Lenstra} and proof is given.

In Theorem \ref{th5} below we give a formula for the smallest number of generators of any ring $M_2(\mathbb{Z})^{k}$. We reduce this problem to the one about the smallest number of generators of the algebras $M_2(\mathbb{F}_p)^k$, where $\mathbb{F}_p$
denotes the field of $p$ elements, $p$ being prime. Let $q$ be a power of a prime; we show in Theorem \ref{th4} that
\begin{equation}\label{eeq1}
\text{gen}_{m,2} \left( \mathbb{F}_q \right) =
\frac{q^{4m-1}+q^{2m} -q^{3m}-q^{3m-1}}{q^2-1}.
\end{equation}
Theorem \ref{th5} is proved by showing that $\text{gen}_{m,2}\left( \mathbb{Z} \right) = \text{gen}_{m,2}\left( \mathbb{F}_2 \right)$ for all $m \geq 3$ which together with (\ref{eeq1})
gives (\ref{salvo}). We do not know how to solve the following problems.

\begin{problem}\label{q1} Find a formula for $\gen_{m,n}(\mathbb{Z})$ for $n \geq 3$.
\end{problem}

\begin{problem}\label{prob2} Find a formula for
$\gen_{m,n} \left( \mathbb{F}_q \right)$ for all $m,n \geq 2$ and $q$ a positive power
of a prime number.
\end{problem}

\begin{problem}\label{prob3}
Is it true that $\gen_{m,n} \left( \mathbb{Z} \right) = \gen_{m,n} \left( \mathbb{F}_2 \right)$
for all $m,n \geq 2$?
\end{problem}


We compute $\gen_{m,n}(\mathbb{F}_2)$ in this paper for $n=2$ and our computation gives a strategy how to deal with all other $n \geq 3$.

When $n \geq 3$, our results are less precise than (\ref{eeq1}). Namely, in Theorem \ref{th2} we prove the following asymptotic formula: let $m,n \geq 2$ be fixed and
$q \to \infty$, then
$\text{gen}_{m,n} \left( \mathbb{F}_q \right)$ is strictly bounded above by and is asymptotically equivalent to
\begin{equation}\label{eeq3}
\left( q-1 \right) q^{(m-1)n^2} \prod_{k=1}^{n} \left( 1 - q^{-k} \right)^{-1} \text{.}
\end{equation}
This number is bounded above by $3.463 \, \left( q-1 \right) q^{(m-1)n^2}$.
It follows that
\[
\text{gen}_{m,n}(\mathbb{Z}) \leq \text{gen}_{m,n} \left( \mathbb{F}_2 \right) <
3.463 \times 2^{(m-1)n^2}  \text{.}\]

On comparing the exact formula (\ref{eeq1}) with the asymptotic formula (\ref{eeq3}) when $n = 2$,
we see that the latter holds for $\gen_{m,2}(q)$ as $m,q \to \infty$; thus in this case (\ref{eeq3}) is true under wider conditions than $m$ fixed as $q \to \infty$. This phenomenon may occur in other cases.

Theorem \ref{gap} states that for a fixed $n \geq 2$ and commutative ring $R$, if the smallest number $r_{m,n}(R)$ of generators of the $R$-algebra $M_n(R)^{m}$ is less than $r_{m+1,n}(R)$, then
$r_{m+1,n}(R) - r_{m,n}(R) = 1$. In addition, Theorem \ref{gap} provides the following non-optimal lower bound:
\begin{equation*}
\gen_{m+1,n}(R) \geq 2 \, \gen_{m,n}(R) \text{.}
\end{equation*}

\paragraph*{Acknowledgements.} The first author was supported by NSF Grant DMS-0456185.
The second author thanks Said Sidki for introducing him to the Hall's paper \cite{Hall}, for many stimulating conversations, without which this paper could not have been conceived. The second author thanks Nigel Boston for stimulating conversations, in particular for his remarks that contributed to the discovery of formulas (\ref{eeq1}), (\ref{eeq3}), and Table \ref{tablegen} below, for performing a crucial MAGMA computation \cite{Boston}, and for his warm hospitality during the second author's visit to the University of Wisconsin-Madison. The second author thanks Derek Robinson for stimulating conversations and his warm hospitality during the second author's visit to the University of Illinois at Urbana-Champaign. The second author thanks John Tate for a useful conversation that contributed to the discovery of Table \ref{tablegen} below. We thank Martin Kassabov for finding errors in formulas (\ref{eeq1}) and (\ref{eeq3}) in a previous version of this paper. We thank Marcin Mazur for improving Theorem \ref{modernch} below. The second author enjoyed his wonderful hospitality during his visit at SUNY Binghamton. We thank Hendrik Lenstra for generously sharing his ideas with us and especially for communicating to us the crucial Theorem \ref{Lenstra1}. We thank Christopher Hillar for finding a typo the in the proof of Theorem \ref{Lenstra1}.

\section{Preliminary Results}

\subsection{The structure of generators for finite direct sums of matrix algebras over a field}

\begin{lemma}\label{lem1}
Let $F$ be a field, and $k\geq 2$ and $n \geq 1$ be integers. Then any $k$ pairs of matrices
$(A_1, B_1), \ldots (A_k, B_k) $ generate the algebra $M_n(F)^2$ if and only if the following
two conditions are
satisfied:
\begin{enumerate}
\item\label{it1} $A_1, \ldots A_k$ generate $M_n(F)$ as an $F$-algebra, and the same
is true of $B_1, \ldots B_k$.
\item\label{it2} There is no $C \in GL_n(F) $ such that $A_i = C^{-1} B_i C$ for $i = 1, \ldots, k$.
\end{enumerate}
\end{lemma}

\begin{proof}
We only need to prove that these conditions are sufficient. Consider the $F$-subalgebra
$\mathcal{A}$ of $M_n(F)^2$ generated by $(A_1,B_1), \ldots, (A_k,B_k)$. Suppose that
\[\mathcal{A} \subsetneqq M_n(F)^2.\]
We will show that each $A_i$ is conjugate to $B_i$
by the same matrix with entries in $F$. Consider
the projections
\[pr_1, \,pr_2: \mathcal{A} \to M_n(F)\]
to the first and second components of $M_n(F)^2$, respectively.
Then
\[pr_1(\mathcal{A}) = pr_2(\mathcal{A}) = M_n(F)\]
because of Condition
\ref{it1} of the lemma. Let
$\mathcal{I}_i = \ker (pr_i) $, $i = 1,2$.

We claim that $\mathcal{I}_1 \neq \mathcal{I}_2$. Otherwise,
\[
\mathcal{I}_1 = \mathcal{I}_2 = \mathcal{I}_1 \cap \mathcal{I}_2 = \{0\}.\]
Therefore,
in addition to being onto, each of $pr_1$ and $pr_2$ is an embedding, and hence
an algebra isomorphism. We may turn $F^n$ into a simple left $\mathcal{A}$-module in two ways according to whether $\mathcal{A}$ acts on $F^n$ via $pr_1$ or $pr_2$. These two modules must be isomorphic.
Let $C$ be the matrix corresponding to this isomorphism in the standard basis of $F^n$. Then $C$ is invertible, and for all $a \in \mathcal{A}$ and $v \in F^n$, we have
\[
C \, pr_1(a) v = pr_2(a) Cv.
\]
Let $a_i = (A_i,B_i)$ where $i = 1, \ldots n$. Then $CA_iv = B_iCv$.
We conclude that $CA_i = B_iC$, which contradicts Condition \ref{it2} of the lemma.

Therefore, $\mathcal{I}_1  \neq  \mathcal{I}_2$, and since $\mathcal{I}_1$ and $\mathcal{I}_2$ are maximal ideals of $\mathcal{A}$, we conclude that $\mathcal{I}_1  +  \mathcal{I}_2 = \mathcal{A}$. Hence,
by the Chinese Remainder Theorem, we have an algebra isomorphism
\[
\mathcal{A} / \mathcal{I}_1 \cap \mathcal{I}_1 \cong \mathcal{A} / \mathcal{I}_1 \oplus
\mathcal{A} / \mathcal{I}_2 = M_n(F)^2.
\]
Then the $F$-dimension counting shows that
$\mathcal{A} = M_n(F)^2$. This final contradiction proves the lemma.
\end{proof}

We generalize this lemma in the next theorem.

\begin{theorem}\label{th1}
Let $F$ be a field, and $k,m\geq 2$ and $n \geq 1$ be integers. Then any $k$ (double-subscripted) $m$-tuples of matrices
$(A_{11}, \ldots,  A_{1m}), \ldots  (A_{k1}, \ldots,  A_{km})$ generate the algebra $M_n(F)^m$ if and only if the following
two conditions are
satisfied:
\begin{enumerate}
\item\label{it3} For any $i = 1, \ldots, m$, the matrices $A_{1i}, \ldots, A_{ki}$ generate $M_n(F)$ as an $F$-algebra.
\item\label{it4} There is no $i \neq j$ such that there is $C \in GL_n(F) $ such that
\[A_{1i} = C^{-1} A_{1j} C, \ldots,  A_{ki} = C^{-1} A_{kj}C.\]
\end{enumerate}
\end{theorem}

\begin{proof} We only need to prove that the conditions of the theorem are sufficient. Suppose the theorem is false, and let $a_j = (A_{j1}, \ldots,  A_{jm})$ where $j = 1, \ldots, k$ provide a counterexample. Let $\mathcal{A}$ be the $F$-subalgebra of $M_n(F)^m$ generated by
$a_1, \ldots, a_k$. Let $pr_i : \mathcal{A} \to M_n(F)$, $i = 1, \ldots, m$, be the projection map onto the
$i$th component, and let $\mathcal{I}_i = \ker (pr_i)$. All the ideals $\mathcal{I}_i$ are maximal
in $\mathcal{A}$. Therefore,  they may not be all different, because otherwise by the Chinese Remainder Theorem,
$\mathcal{A} = M_n(F)^n$, a contradiction.
We may assume that
$\mathcal{I}_1 = \mathcal{I}_2$ without loss of generality. Therefore,
\[
\mathcal{I}_1 = \mathcal{I}_2 = \mathcal{I}_1 \cap \mathcal{I}_2 =
\]
\[
\left\{ \left(0,0,M_1, \ldots, M_{m-2} \right) \mid M_1, \ldots, M_{m-2} \in M_n(F)\right\}.
\]
Let
\[
pr_{12} : M_n(F)^m \to M_n(F)^2,~(D_1, \ldots, D_m) \mapsto (D_1, D_2).
\]
Then by Lemma
\ref{lem1},
\[
\mathcal{A}' = pr_{12} (\mathcal{A}) = M_n(F)^2.
\]
At the same
time,
\[
\mathcal{I}_i' = pr_{12} (\mathcal{I}_i) = \{ 0\} ~ \text{for} ~ i = 1,2.
\]
Therefore,
$\mathcal{I}_1'$ is neither a maximal ideal of $\mathcal{A}'$ nor it is $\mathcal{A}'$,  contradicting that
$\mathcal{I}_1$ is a maximal ideal of $\mathcal{A}$.
\end{proof}

Condition 1 in Theorem \ref{th1} for $n \geq 2$ may be verified by the Theorem of Burnside from \cite{Burnside}.
In today's terminology, it states that a collection $S$ of matrices generates the matrix algebra $M_n(F)$ over a field $F$ if and only if the matrices in $S$ do not have a common eigenspace over an algebraic closure of $F$.
Other than using the definition, we do not know how to decide whether two pairs of matrices
are conjugate by the same matrix.

\subsection{The irrelevance of the identity matrix in any generating set of a matrix algebra}

The next result shows that the identity element may always be removed from
any generating set of a finite direct sum of integer matrix rings. We prove a more general result.
For the convenience of the reader, we preface the proof by stating the following two standard facts.
\begin{facts}\label{naive}
1. Let $R$ be a ring, and let $A_1, \ldots, A_k$ be $R$-algebras. Then any two-sided ideal of
$\bigoplus_{i=1}^k A_i$ is of the form $\bigoplus_{i=1}^k J_i$ for some
two-sided ideals $J_i$ of $A_i$.

2. Any two-sided ideal $\mathcal{J}$ of any matrix algebra
$M_n(R)$ is of the form $M_n(J)$ where $J$ is the two-sided ideal of $R$
generated by all the entries of all matrices
in $\mathcal{J}$.
\end{facts}

\begin{theorem}\label{carefree}
Let $R$ be a commutative ring, and
$A$ a finite direct sum of matrix algebras with entries in $R$. Then
\begin{enumerate}
\item\label{1it} A set $S$
generates $A$ as an $R$-algebra if and only if the set
$ S \cup \left\{ 1_A \right\} $ does.
\item\label{2it} Any maximal subalgebra of $A$ contains the identity element.
\end{enumerate}
\end{theorem}

This theorem is false if we allow $1$-by-$1$ summands in the direct sum because, for example, for a field $F$ the algebra $F \oplus F$ contains maximal subalgebras not containing the element $(1,1)$.

\begin{proof}[Proof of Theorem \ref{carefree}] \ref{1it}. Suppose that Part \ref{1it} of the theorem is false. Then the set $S \cup \left\{1_{A}\right\}$
generates $A$, while $S$ does not. Let $\langle S \rangle$ denote the $R$-subalgebra of $A$
generated by $S$. Then
\begin{equation}\label{Rostik}
A = 1_{A} R + \langle S \rangle \text{.}
\end{equation}
Furthermore, (\ref{Rostik}) implies that $\langle S \rangle$ is a two-sided ideal of $A$.
Therefore (\ref{Rostik}) implies that $A / \langle S \rangle$ is an epimorphic image of
$R$; hence $A$ is a commutative ring.
On the other hand, Facts \ref{naive} tell us that $A / \langle S \rangle$
is isomorphic to a finite direct sum of matrix rings with entries in some nonzero commutative rings.
Therefore, $A$ is not commutative, a contradiction.

\ref{2it}. Let $\mathcal{A}$ be a maximal subalgebra of $A$, and let $\mathcal{B}$
be the subalgebra of $A$ generated by $\mathcal{A} \cup \{ 1_{A} \}$. Then
$\mathcal{B} \subsetneqq A$ by Part \ref{1it}. Since in addition $\mathcal{A}$ is a maximal subalgebra contained in
$\mathcal{B}$, we conclude that $\mathcal{A} = \mathcal{B}$.
\end{proof}

The following result has been improved by Marcin Mazur \cite{Mazur} who showed that the polynomial
in it may be elegantly expressed as a determinant of a commutator.

\begin{theorem}\label{modernch}
Let $R$ be a commutative ring with the group of invertible elements $U(R)$.
The matrices $A,B \in M_2(R)$ generate $M_2(R)$ as an $R$-algebra if and only if
\begin{equation}\label{wonderdet}
\det
\left(
\begin{array}{cccc}
1 & 0 & 0 & 1\\
a_{11} & a_{12} & a_{21} & a_{22}\\
 b_{11} & b_{12} & b_{21} & b_{22} \\
&\text{Flatten}(AB)\\
\end{array}
\right) = \det(AB - BA) \in U(R)  \, \text{,}
\end{equation}
where $\text{Flatten}(AB)$ is the matrix $AB$ written as a row-vector with $4$ components,
similarly to the first $3$ rows that are $\text{Flatten}\left(I_2 \right)$,
$\text{Flatten}(A)$, and $\text{Flatten}(B)$, respectively.

The polynomial $\det(AB - BA)$ is irreducible over any field.
\end{theorem}

\begin{proof} By Part \ref{1it} of Theorem \ref{carefree}, it suffices to prove that (\ref{wonderdet}) holds
for any $A,B \in M_2(R)$ such that $I_2, A,B$ generate the algebra $M_2(R)$.
Then there exist some words $w_1, \ldots, w_k$ in $A$ and $B$ such that
$M_2(R) = R\, I_2 + \sum_{i=1}^k  R w_i$. By Cayley-Hamilton
theorem, we may assume that each $w_i$ contains no powers of $A,B$ higher than the $1$st.
In addition, because of
\[
BA = (A+B)^2 - A^2 - B^2-AB=
\]
\[
\tr(A+B) \,(A+B) - \det(A+B) \, I_2 - \tr A  \,A + \det A  \,I_2 - \tr B \, B + \det B  \,I_2 -AB \text{,}
\]
we may assume that
$\left\{ w_1, \ldots, w_k \right\} \subseteq \left\{I_2, A,B, AB\right\}$.
Then
\begin{equation}\label{Said}
M_2(R) = R\, I_2 + R\,A + R\,B + R \, AB \text{.}
\end{equation}
Let $f(A,B) = $
\[
f \left(  a_{11},a_{12}, a_{21},a_{22},b_{11}, b_{12}, b_{21},b_{22} \right) =
\det
\left(
\begin{array}{cccc}
1 & 0 & 0 & 1\\
a_{11} & a_{12} & a_{21} & a_{22}\\
 b_{11} & b_{12} & b_{21} & b_{22} \\
&\text{Flatten}(AB)\\
\end{array}
\right).
\]
Then, (\ref{Said}) and $4 = \dim_{R} M_2(R) =  \# \left\{  I_2, A,B,AB \right\}$
imply that $A,B \in M_2(R)$ generate the $M_2(R)$ as an $R$-algebra if and only if
$f(A,B)$ is an invertible element of the ring $R$.

\medskip

Finally we will show that the two polynomials in (\ref{wonderdet}) are the same polynomial, and this polynomial is irreducible.

To prove irreducibility of the polynomial $f$, we calculate the following $2$ specializations:
\begin{enumerate}
\item
$f \left( a_{{11}},0,1,0,0,1,b_{{21}},0 \right) =
 -1+b_{{21}}{a_{{11}}}^{2}$.
\item
$f \left( 0,a_{{12}},1,0,0,1,b_{{21}},0 \right) =
- \left( a_{{12}}b_{{21}}-1 \right) ^{2}$.
\end{enumerate}
The first specialization is
an irreducible cubic, and the second one is the square of an irreducible quadratic.
These claims are easily established by the method of undetermined coefficients.

Now we will show that $f(A,B) = \det(AB - BA)  = \det[A,B]$. Let $R = \mathbb{C}$. Suppose that the matrices $A,B \in M_2(\mathbb{C})$ do not generate the algebra
$M_2(\mathbb{C})$. Then
\begin{equation}\label{var1}
f(A,B) = 0\text{.}
\end{equation}
Then by Theorem \ref{carefree}, the matrices $I_2, A,B$ do not generate the the algebra $M_2(\mathbb{C})$ as well. Then the subalgebra $S$ that they generate is at most $3$-dimensional, and in particular $\dim_{\mathbb{C}} S/\text{rad}(S) \leq 3$. Therefore, by Wedderburn theory, the algebra $S/\text{rad}(S)$ is commutative. In other words, $[A,B] \in \text{rad}(S)$.
Every element of $\text{rad}(S)$ is nilpotent, and the square of every nilpotent matrix in $M_2(\mathbb{C})$ is zero. Therefore $[A,B]^2 = 0$. By Cayley-Hamilton theorem,
$0 = [A,B]^2 - \text{tr}[A,B]\,[A,B] + \det[A,B] I_2 = \det[A,B] \, I_2$, so that
\begin{equation}\label{var2}
\det[A,B] = 0\text{.}
\end{equation}
Therefore, by (\ref{var1}) and (\ref{var2}), the set of zeros of $f(A,B)$ is contained in the set of zeros of $\det[A,B]$. Both polynomials have the same degree, and the polynomial $f$ is irreducible. Therefore, by Hilbert's Nullstellensatz
\[f(A,B) = c \det[A,B]
\]
for some $c \in \mathbb{C}$. By evaluating both polynomials at $A = E_{12}$ and $B = E_{21}$ we find that $c = 1$. Therefore, the polynomials $f(A,B)$ and $\det[A,B]$ are equal over $\mathbb{Z}$ and therefore over any commutative ring.
\end{proof}

We remark that Theorem \ref{modernch} and Burnside's Theorem \cite{Burnside} imply that the equation $\det[A,B] = 0$ describes the set of pairs $(A,B)$ of $2$-by-$2$ matrices with entries in a given field $F$ such that $A$ and $B$ have a common eigenvector over an algebraic closure of $F$.

The polynomial $\det[A,B]$ has some interesting properties. For example, over any commutative ring $R$
we have
$\det[A,B] = \det[A - s\,I_2,\, B- t\,I_2]$ for all $s,t \in R$.

\bigskip

The question of whether the given $A,B \in M_n(\mathbb{Z})$ generate $M_n(\mathbb{Z})$
as a ring has an algorithmic solution (see \cite[Theorem 2.9]{PetSid}). It consists of deciding whether the rows of a certain nonsquare integer matrix span the full lattice. An application
of Theorem \ref{carefree} has helped us simplify this algorithm to (\ref{wonderdet})
in the case of $2$-by-$2$ matrices.

\subsection{The smallest number of generators for finite direct sums of matrix algebras of different sizes}
The result of this subsection states that the smallest number of generators of a finite direct sum of matrix algebras of different sizes is the maximum of the smallest number of generators in all the sums of the terms of the same size.

\begin{theorem}\label{th6} Let $R$ be a commutative ring. Let
$n_1, \ldots, n_k \geq 1$ be pairwise different positive integers, and
$m_1, \ldots, m_k$ positive integers. Suppose that for any $i = 1, \ldots, k$,
the $R$-algebra $M_{n_i}(R)^{m_i}$ has $s_i$ generators. Then the $R$-algebra
$\bigoplus_{i = 1}^k M_{n_i}(R)^{m_i}$ has $\max \{ s_1, \ldots, s_k\}$ generators.
\end{theorem}

\begin{proof}
We may assume that $s = s_1 \geq s_2 \geq \ldots \geq s_k$. Consider
the free associative $R$-algebra
$\mathcal{R} = R \{ x_1, \ldots, x_s\}$ in free $s$ variables. Since
the algebra $M_{n_i}(R)^{m_i}$ has $s_i$ generators, there exists an algebra
epimorphism $\pi'_i: R \{ x_1, \ldots, x_{s_i}\} \to M_{n_i}(R)^{m_i}$ such that
$\pi'_i \left( x_1 \right), \ldots, \pi'_i \left( x_{s_i} \right)$ generate $M_{n_i}(R)^{m_i}$.
Let $\pi_i: \mathcal{R} \to M_{n_i}(R)^{m_i}$ be obtained by composing
$\pi'_i$ with the evaluation map
$\mathcal{R}\{ x_1, \ldots, x_{s} \} \to R \{ x_1, \ldots, x_{s_i} \}$ given by
$f\left(x_1, \ldots, x_s\right) \mapsto f\left(x_1, \ldots, x_{s_i}, 0, \ldots,0\right)$.
Let $J_i = \ker \pi_i$.

We note that these ideals are pairwise not contained in each other.
To see that, it is enough to show that $J_1 \nsubseteq J_2$ and $J_2 \nsubseteq J_1$. If this
were not true, then for example, $J_1 \subseteq J_2$, so that the $R$-algebra
$M_{n_2}(R)^{m_2}$ is an epimorphic image of $M_{n_1}(R)^{m_1}$. Let $M$ be a maximal ideal of $R$,
and $F$ the field $R/M$. Then then the $F$-algebra
$M_{n_{_2}}(F)^{m_{_2}} \cong M_{n_{_2}}(R)^{m_{_2}} \otimes_{R} F$ is an epimorphic image of
$M_{n_{_1}}(F)^{m_{_1}} \cong M_{n_{_1}}(R)^{m_{_1}} \otimes_{R} F$. This is impossible because $n_1 \neq n_2$ and because of Facts \ref{naive}.

In view of the Chinese Remainder Theorem, if we show that
$J_i + J_j = \mathcal{R}$ whenever $i \neq j$, then
\[
\frac{\mathcal{R}}{\bigcap_{i=1}^k J_i} \cong \bigoplus_{i=1}^k M_{n_i}(R)^{m_i}\text{,}
\]
and therefore the latter
$R$-algebra also has $s$ generators
\[
x_1 + \bigcap_{i=1}^k J_i, \ldots, x_s + \bigcap_{i=1}^k J_i \text{.}
\]
The element $1 + \bigcap_{i=1}^k J_i$ is not needed in this generating set by Part \ref{1it} of Theorem \ref{carefree}.

It suffices to show that $J_1 + J_2 = \mathcal{R}$. If $J_1 + J_2 \subsetneqq \mathcal{R}$ then
for $j = 1,2$
\begin{equation}\label{eq5}
M_{n_j}(R)^{m_j} \cong \mathcal{R} / J_j \stackrel{\text{onto}}{\longrightarrow}
\mathcal{R} / (J_1 + J_2) \cong M_{n_j}(R/M_j)^{m_j}  = A_j
\end{equation}
for some proper ideals $M_1, M_2$ of $R$. {\em However, the $R$-algebras $A_1$ and $A_2$
are not isomorphic, in contradiction with (\ref{eq5})}. Indeed,
if $M = M_1 = M_2$, then
$A_1 \cong A_2$ implies that $M_{n_1}(R/M)$ is isomorphic to a finite direct
sum of copies of $M_{n_2}(R/M)$. Since $n_1 \neq n_2$, this leads to a contradiction as before,
by tensoring over $R/M$ both sides of this purported direct decomposition with
any residue field of $R/M$. Finally, if $M_1 \neq M_2$, then the $R$-annihilators of
$A_1$ and $A_2$ are different.
\end{proof}

\subsection{The size of gaps in the sequence $\gen_{m,n}(R)$ of the smallest number of generators is at most $1$}
\begin{definition} Let $R$ be a commutative ring, and $m,n \geq 1$ integers. We define the set
$G_{m,n}(R) \subseteq M_n(R)^m$ to consist of the elements
$(A_1, \ldots, A_m)$ such that $A_1, \ldots, A_m$ generate $M_n(R)$ as an $R$-algebra.

We introduce the sequence
$ \{ \gen_{m,n}(R) \}$ defined by the property that $\gen_{m,n}(R)$ is the integer such that
the $R$-algebra $M_n(R)^{\gen_{m,n}(R)}$ admits $m $ generators, while $M_n(R)^{1+\gen_{m,n}(R)}$
needs at least $m+1$ generators. If any $R$-algebra $M_n(R)^k$ admits $\leq m$ generators,
then put $\gen_{m,n}(R) = \infty$.
\end{definition}


We see that for any ring $R$, we have $\gen_{m,n}(\mathbb{Z}) \leq \gen_{m,n}(R)$. The set $G_n(R) = G_{2,n}(R)$ for $n \geq 2$ has been introduced in \cite{PetSid}.
This set is nonempty for any ring because $M_n(R)$ as an $R$-algebra has $2$ generators, for example,
$X = E_{1n} + \sum_{i=1}^{n-1} E_{i+1,i}$ and $Y = E_{11}$,
since $E_{ij} = X^{i-1} Y X^{n-j+1}$.

\begin{theorem}\label{gap}
As $l \to \infty$, whenever the smallest number of generators of $M_n(R)^l$ increases, the increment is exactly $1$. Furthermore, for any $m,n \geq 2$
\begin{equation}\label{double}
\gen_{m+1,n}(R) \geq 2 \, \gen_{m,n}(R) \text{.}
\end{equation}
\end{theorem}

\begin{proof}
With the help of the matrices $X = E_{1n} + \sum_{i=1}^{n-1} E_{i+1,i}$ and $Y = E_{11}$, we see that if the algebra $M_n(R)^l$ has $k$ generators
\[
a_1 = (A_{11}, \ldots, A_{1l}), \ldots, a_k = (A_{k1}, \ldots, A_{kl}) \text{,}
\]
then $M_n(R)^{l+1}$ has the following $k+1$ generators:
\[
a_1' = (A_{11}, \ldots, A_{1l}, 0), \ldots, a_{k-1}' = (A_{k-1,1}, \ldots, A_{k-1,l}, 0),
\]
\[
a_k' =  (A_{k1}, \ldots, A_{kl}, Y),~a_{k+1}' = (0, \ldots, 0, X).
\]
The last claim is true because the last components
of $a_k' a_{k+1}'$ and $a_{k+1}'$ generate the $R$-algebra $M_n(R)$, and in addition,
$\left(a_{k+1}' \right)^{n} a_k'= (0, \ldots, 0, Y)$.

Next, we construct $k+1$ generators $a_1'', \ldots, a_{k+1}''$ for the $R$-algebra $M_n(R)^{2l}$
as follows: $a_1'', \ldots, a_k''$ are respectively the juxtaposed $a_1$ with $a_1$, $\ldots$, $a_k$ with $a_k$, i.e.~
\[
a_1'' = \left( A_{11}, \ldots, A_{1l}, A_{11}, \ldots, A_{1l} \right)\text{,} \ldots\text{,}~
a_k'' = \left( A_{k1}, \ldots, A_{kl}, A_{k1}, \ldots, A_{kl} \right) \text{.}
\]
Finally, each the first $l$ components $a_{k+1}$ is zero, and each of the remaining
$l$ components is $I_n$, i.e.
\[
a''_{k+1} =(
\underbrace{0, \ldots,  0}_{l~\text{components}},~
\underbrace{I_n, \ldots,  I_n}_{l~\text{components}}
)
\]
We see that the $R$-linear combinations of products of $a''_1, \ldots, a''_{k+1}$,
span $M_n(R)^{2l}$.
\end{proof}

The bound (\ref{double}) is general. However, it is probably never sharp.

\subsection{Two local-global results}
The two results below will be later applied to finite direct sums of integer matrix rings.

Firstly, we generalize Theorem 2.4 of \cite{PetSid}.
\footnote{It was originally communicated to the second author by David Saltman (University of Texas at Austin).}

\begin{theorem}\label{th3}
Let $R$ be a ring whose additive group $(R,+)$ is finitely generated. Then a subset $S$ generates $R$ as a ring if and only if the reduction of $S$ modulo every prime $p$ generates the ring $R/pR$. If $(R,+)$
is infinite, then in general no prime number may be omitted in the previous sentence.
\end{theorem}

\begin{proof} We see that if $H$ is an additive subgroup of $R$,
then $H = R$ if and only if the reduction of $H$ modulo any prime $p$ is $R/pR$.
The first statement in the theorem is proved by applying this fact to the additive group of the subring
generated by $S$.

Let $(R,+)$ be infinite, and suppose that we need not consider the reduction modulo some prime $p_0$. Let $S$ be a generating set for
$R$ as a ring. Then the set $p_0 S$ clearly does not generate the ring $R$. At the same time, modulo every prime $p \neq p_0$ the set $S$ generates
the ring $R/pR$.
\end{proof}

\underline{Therefore,} Theorems \ref{th1} and \ref{th3} imply that to prove that a set
\[
S = \left\{ s_1 = \left( A_{11}, \ldots, A_{1n} \right),
s_2 = \left( A_{21}, \ldots, A_{2n} \right),
\ldots ,  s_m = \left( A_{m1}, \ldots, A_{mn} \right)
\right\}
\]
generates the ring $M_n(\mathbb{Z})^m$ it is necessary and sufficient to prove that
\begin{itemize}
\item Each of the $n$ sets $\left\{ A_{11}, A_{21}, \ldots, A_{m1}  \right\}$,
$\left\{ A_{12}, A_{22}, \ldots, A_{m2}  \right\}$, $\ldots$, \\
$\left\{ A_{1n}, A_{2n}, \ldots, A_{mn}  \right\}$ generates the ring
$M_n(\mathbb{Z})$.
\item For any prime $p$, any two of the $n$ ordered $m$-tuples
$\left( A_{11}, A_{21}, \ldots, A_{m1}  \right)$, \\
$\left( A_{12}, A_{22}, \ldots, A_{m2}  \right)$, $\ldots$,
$\left( A_{1n}, A_{2n}, \ldots, A_{mn}  \right)$ are not conjugate to each other modulo $p$ by the same matrix.
\end{itemize}

\bigskip

Secondly, we reproduce the proof of theorem of Lenstra from \cite{Lenstra_C} 
with some minor stylistic changes.
\begin{theorem}\label{Lenstra}
Let  $R$  be a ring whose additive group is finitely generated.
For each prime $p$,  let
            \[r(p)\]
be the smallest number of generators of
$R/pR$ as a ring, and let
            \[r(0)\]
be the smallest number of generators of
$ R\otimes_{\mathbb{Z}} \mathbb{Q}$ as a $\mathbb{Q}$-algebra. Finally, let
            \[r\]
be the smallest number of generators of  $R$  as a ring. Then
\begin{enumerate}
\item\label{it1} For each  prime number  $p$  we have $r(0) \le r(p) \le r$.
\item\label{it2} For all but finitely many prime numbers we have $r(0) = r(p)$.
\item\label{it3} If there is a prime number  $p$  such that  $r(0) < r(p)$,  then
                        \[r = \max\{r(p): p~ \text{prime}\}.\]
\item\label{it4} If $r(0) = r(p)$  for all  $p$,  then
            either  $r = r(0)$  or  $r = r(0) + 1$.
\end{enumerate}
\end{theorem}

In particular since $(1,1,0)$ and $(0,1,1)$ generate the rings $\mathbb{Z}^3$ and $\mathbb{F}_p^3$
for all primes $p$, and the $\mathbb{Q}$-algebra $\mathbb{Q}^3$, and because of (\ref{cute}) below,
we conclude that
\begin{equation}\label{cute1}
\gen_{m,1}(\mathbb{Z}) = \gen_{m,1}(2) = 2^m-1.
\end{equation}
A similar strategy will be used to find $\gen_{m,2}(\mathbb{Z})$, but the computations will be longer.

\begin{proof}[Proof of Theorem \ref{Lenstra}]
If a set $S$ generates the ring $R$, then  for every prime $p$ the reduction of $S$ generates the ring  $R/pR$. Let $p$ be a prime, and $T$ a generating set of the ring $R/pR$. Then the preimage of $T$ under the reduction map $R \to R/pR$ generates a subring $R_0$ of $R$ of $R$ of finite index as an additive subgroup. Therefore $R_0 \otimes 1$ generates $R\otimes_{\mathbb{Z}} \mathbb{Q}$ as a $\mathbb{Q}$ algebra. This proves Part \ref{it1}.

In view of Part \ref{it1}, to prove Part \ref{it2} we need to show that $r(0) \geq r(p)$ for all but finitely many prime numbers $p$. Let $a_1, \ldots, a_t$ be a generating set of the $\mathbb{Q}$-algebra
$R\otimes_{\mathbb{Z}} \mathbb{Q}$. Then for a sufficiently large integer $s_i$ we may write
$s_i a_i = \sum_{j=1}^k b_j \otimes \alpha_{ij}$ for some $b_j \in R$ and integers $\alpha_{ij}$. Then
the $t$ elements $c_i = \sum_{j=1}^k \alpha_{ij} b_j$ generate the subring $R_1$ of $R$ of finite index $t_1$. Then for any prime $p$ not dividing $t_1$, the image of $R_1$ under the reduction map
$R \to R/pR$ is $R/pR$. This proves Part \ref{it2}.

Parts \ref{it3} and \ref{it4} will be proved simultaneously, and this constitutes the most important part of the theorem. Let
\begin{equation}\label{bound}
m= r(0),~\, n = \max \left\{m+1~\text{and all}~ r(p)~\text{where}~ p~ \text{is prime}  \right\}.
\end{equation}
Starting from $m$ elements $x_1, \ldots, x_m$ of $R$ that give $\mathbb{Q}$-algebra generators of  $R\otimes_{\mathbb{Z}} \mathbb{Q}$, we will construct $n$ generators of the ring $R$. The elements $x_1, \ldots, x_m$ will be replaced one by one by elements $y_1, \ldots, y_m$ of $R$, to be constructed carefully by induction. After that we will construct the elements
$x_{m+1}, \ldots, x_{n}$ such that the elements $y_1, \ldots, y_m, x_{m+1}, \ldots, x_{n}$
generate the ring $R$.

Let $B_0$ be the set of primes $p$ such that the ring $R/pR$ does not admit $n-1$ generators. This set is finite by the definition of $n$ in (\ref{bound}) and by Part \ref{it2}. At the $k$th stage, we have elements $y_1, \ldots, y_k$ of $R$ that together with $x_{k+1}, \ldots, x_m$ give $\mathbb{Q}$-algebra generators of $R\otimes_{\mathbb{Z}} \mathbb{Q}$. The finite set $B_k$ is defined as the set of all primes $p$ such that
the elements $y_1, \ldots, y_k$ when taken modulo $p$ \underline{do not form} part of a set of $n-1$ generators for the ring $R/pR$. It will also be true that the the $k$th stage, for every prime $p$, the elements
$y_1, \ldots, y_{k}$ taken modulo $p$ \underline{form} part of a set of $n$ ring generators for $R/pR$. (This is also correct for $k = 0$ because by (\ref{bound}), we have $n \geq r(p)$ for any prime $p$.)

Next we explain how the argument passes from stage $k$ to stage $k+1$ to construct $y_{k+1}$. Firstly, we construct $z \in R$ as follows. For any
$p \in B_k$ we may choose $z_p \in R$ such that $y_1, \ldots, y_k, z_p$ taken modulo $p$ form part of a set of $n$ ring generators for $R/pR$. Let
\[
z = \sum_{p \in B_k} z_p \, \prod_{p \neq p' \in B_k}p'.
\]
We see that $z \equiv z_p ~\mod\,p$ for every $p \in B_k$. Hence for \underline{every} prime $p$, the elements $y_1, \ldots, y_k, z$ taken modulo $p$ form part of a set of $n$ ring generators for $R/pR$. Therefore, by Theorem \ref{th3}, the elements
$y_1, \ldots, y_k, z$ form part of a set of $n$ ring generators for $R$.

The algebra $R\otimes_{\mathbb{Z}} \mathbb{Q}$ is isomorphic as a $\mathbb{Q}$-vector space to some
$\mathbb{Q}^s$. The Euclidean norm on $\mathbb{Q}^s$ induces a norm on $R\otimes_{\mathbb{Z}} \mathbb{Q}$ making it a linear normed space over a normed field; multiplication in $R\otimes_{\mathbb{Z}} \mathbb{Q}$ is continuous as a bilinear map. Let $P$ be an integer divisible by all the primes in $B_k$ and so large that the element
$x_{k+1} + z/P$ is sufficiently close to $x_{k+1}$ in this norm to ensure that the elements
$y_1, \ldots, y_k, x_{k+1} + z/P, x_{k+2}, \ldots, x_m$ and therefore the elements
$y_1, \ldots, y_k, P\,x_{k+1} + z, x_{k+2}, \ldots, x_m$ generate $R\otimes_{\mathbb{Z}} \mathbb{Q}$
as a $\mathbb{Q}$-algebra. Now let
\[
y_{k+1} = z + P\, x_{k+1}.
\]
We see that the elements $y_1, \ldots, y_{k+1}$ form part of a set of $n$ generators of the ring $R$.

At this point we have constructed $y_1, \ldots, y_m$ which form part of a set of $n$ generators of the ring $R$. If they do not generate the ring $R$, then by Part \ref{it2} and Theorem \ref{th3}, there exists a finite set of primes $B$ such that for any $p \in B$, the elements $y_1, \ldots, y_m$ do not generate the ring $R$. For each $p \in B$, let $x_{m+1,p}\, , \ldots, x_{n,p} \in R$ be such that
the elements $y_1, \ldots, y_k, x_{m+1,p} \, , \ldots, x_{n,p}$ taken modulo $p$ generate the ring
$R/pR$. For $i = m+1, \ldots, n$ let
\[
x_i = \sum_{p \in B} x_{i,p} \, \prod_{p \neq p' \in B}p'.
\]
We have constructed the elements $y_1, \ldots, y_m, x_{m+1}, \ldots, x_{n}$ of $R$ generating $R$ as a ring, and this concludes the proof.
\end{proof}

\section{Applications to generators of finite direct sums of matrix rings }

\subsection{An asymptotic formula for the case of a finite field}
When the group $GL_n \left(  \mathbb{F}_q \right)$ acts on the set $G_{m,n}(\mathbb{F}_q)$ by
conjugating each component of a given $m$-tuple $a$, the centralizer of $a$ is exactly the
intersection $Int$ of the centralizers of all the matrices in the tuple. It follows that $Int$ is the
set of all the nonsingular scalar matrices, and hence the orbit of every point has exactly
$\# PGL_n \left(  \mathbb{F}_q \right)$ elements.
Therefore, Theorem \ref{th1} yields the following formula for the case of a finite field with $q$ elements $\mathbb{F}_q$:
\begin{equation}\label{eq1}
\gen_{m,n}(q) = \gen_{m,n}(\mathbb{F}_q) = \frac{\# G_{m,n}(\mathbb{F}_q)}{\# PGL_n \left(  \mathbb{F}_q \right)}.
\end{equation}
In particular,
\begin{equation}\label{cute}
\gen_{m,1}(\mathbb{F}_q) = \frac{q^m-1}{q-1}.
\end{equation}

Informally, the numerator of the right-hand side of (\ref{eq1}) is obtained by counting the number of $m$-tuples of matrices from
$M_n(\mathbb{F}_q)$, the components of each tuple written below those of the previous one, such that
\begin{itemize}
\item Each vertical cross-section generates
$M_n(\mathbb{F}_q)$ as an $\mathbb{F}_q$-algebra.
\item No two vertical cross-sections are conjugate by the same matrix from $GL_n(\mathbb{F}_q)$.
\end{itemize}

If $m,n$ are fixed
and $q \to \infty$ then by a straightforward generalization of Theorem 2.19 of Petrenko and Sidki
\cite{PetSid}, it follows that
\begin{equation}\label{eq2}
\lim_{q \to \infty} \frac{\# G_{m,n}(\mathbb{F}_q)}{\# M_n(\mathbb{F}_q)^m} = 1.
\end{equation}
Namely, $G_{m,n}(\mathbb{F}_q)$ is obtained from $M_n(\mathbb{F}_q)$ by removing finitely many hypersurfaces, whose number and degrees depend on $m$ and $n$, but not on $q$.
Alternatively, (\ref{eq2}) is a consequence of Theorem 2.19 of \cite{PetSid}. That theorem states that (\ref{eq2}) is true for $m = 2$, and it remains to note that
\[
G_{2,n}\left(\mathbb{F}_q \right) \times M_{n}\left(\mathbb{F}_q \right)^{m-2} \subseteq
G_{m,n}\left(\mathbb{F}_q \right)\text{.}
\]

We pause to restate (\ref{eq2}) in probabilistic terms.
\begin{corollary}
Let $m,n \geq 2$ be fixed. The probability that $m$ matrices
$M_1, \ldots, M_m \in M_n \left( \mathbb{F}_q \right)$ chosen under the uniform
distribution generate the $\mathbb{F}_q$-algebra $M_n \left( \mathbb{F}_q \right)$
tends to $1$ as $q \to \infty$.
\end{corollary}

We see that (\ref{eq1}) and (\ref{eq2}) and
$\# PGL_n \left(\mathbb{F}_q \right) = (q-1)^{-1} \prod_{i = 0}^{n-1} (q^n - q^i)$
imply
\begin{theorem}\label{th2}
If $m,n \geq 2$ are fixed and $q \to \infty$, then
$\gen_{m,n}(q)$ is strictly bounded above by
and is asymptotically equivalent to
\begin{equation}\label{eq8}
\left( q-1 \right) q^{(m-1)n^2} \prod_{k=1}^{n} \left( 1 - q^{-k} \right)^{-1} \text{.}
\end{equation}
In particular, $\gen_{m,n}(q)$
is asymptotically equivalent to
$q^{(m-1)n^2 +1}$ as $m,n\geq 2$ are fixed and $q \to \infty$.
\end{theorem}

\begin{corollary}\label{cor0}
\begin{equation}
\text{gen}_{m,n}(\mathbb{Z}) \leq
2^{(m-1)n^2} \prod_{k=1}^{n} \left( 1 - 2^{-k} \right)^{-1} < 3.463 \cdot 2^{(m-1)n^2}  \text{.}
\end{equation}
\end{corollary}

\begin{proof} Firstly, substitute $q = 2$ in (\ref{eq8}). Then to estimate
$\prod_{k=1}^{\infty} \left( 1 - 2^{-k} \right)^{-1}$ from above,
we use for $x = 1/2$ Euler's Product Formula expressing the reciprocal of our product
as an alternating series:
\begin{equation}\label{eq10}
\prod_{n=1}^{\infty} \left( 1-x^n \right)= 1 +
\sum_{n = 1}^{\infty}(-1)^n \left( x^{^{\frac{3n^2-n}{2}}} +
x^{^{\frac{3n^2+n}{2}}} \right), ~ \, |x|<1 \text{.}
\end{equation}
This allows us to estimate $\prod_{k=1}^{\infty} \left( 1 - 2^{-k} \right)^{-1}$ from above to at least $10^{-6}$ by
summing on the computer the first $3$ terms in the series on the right hand side of (\ref{eq10}) with
$x = 1/2$.
\end{proof}


The methods above, however, tell us nothing when instead of $q \to \infty$ other asymptotics,
such as $n \to \infty$, are considered. While we do not know the answer,
in the proof of Theorem \ref{th4} below, we show that $\#G_{m,2} \left( \mathbb{F}_q \right)$ is asymptotically equivalent
to the number of elements in the ambient set $M_2\left( \mathbb{F}_q \right)^m$ as $m,q \to \infty$.

In the next subsection we obtain Formula (\ref{eq3}) for $\gen_{m,2}(q)$.

\subsection{An exact formula for the case of $2$-by-$2$ matrices with entries in a finite field }

Below we describe in sufficient detail the intersection of any collection of maximal subalgebras of $M_2(\mathbb{F}_q)$. Such subalgebras are of
two types, and we describe each of them in the two paragraphs below.

\underline{\bf{1. The noncommutative maximal subalgebras of $M_2(\mathbb{F}_q)$.}}

\noindent Let $0 \neq v \in \mathbb{F}_q^2$, and let $\mathcal{A}_v \subseteq M_2(\mathbb{F}_q)$
consist of all matrices having $v$ as an eigenvector. Then
$\dim_{\mathbb{F}_q} \mathcal{A}_v = 3$, and since a scaling of $v$ results in the same
subalgebra, we have a one-to-one correspondence between the collection of $\mathcal{A}_v$s
and the $q+1$ points of the projective line $\mathbb{P}_1(\mathbb{F}_q)$.
{\em Hence there are exactly $q+1$ maximal noncommutative subalgebras of
$M_2\left( \mathbb{F}_q  \right)$.}

Below we describe all possible intersections of these maximal subalgebras.
\begin{itemize}
\item If $u \neq v \in \mathbb{P}_1(\mathbb{F}_q)$, then $\mathcal{A}_{u,v} = \mathcal{A}_{u} \cap A_{v}$ has $\mathbb{F}_q$-dimension $2$.
\item If $\{v,v' \} \neq \{w,w' \}$, then $\mathcal{A}_{v,v'} \neq \mathcal{A}_{w,w'}$. Consequently,
there are exactly $\binom{q+1}{2}=(q+1)q/2$ such subalgebras.
\item For any pairwise different $u,v,w \in \mathbb{P}_1(\mathbb{F}_q)$, we have
$\mathcal{A}_u \cap \mathcal{A}_v \cap \mathcal{A}_w = \mathcal{D}$.
\end{itemize}
Indeed, let $u,v \in \mathbb{P}_1(\mathbb{F}_q)$ be different.
Then $\mathcal{A}_{u,v}$ has a basis consisting of two projection operators
on the lines $u$ and $v$, respectively. This proves the first of the above $3$ claims.

The remaining two claims follow from the following
\begin{observation}\label{prop1} Let $F$ be a field, and $a_1,a_2,a_3 \in F^2$ be such that
any two of them are linearly independent.
Then any linear operator that scales $a_1,a_2,a_3$ is a scalar operator.
\end{observation}
\begin{proof} .
Let $a_1$ and $a_2$ be linearly independent, and let $f$ be a linear operator such that $f(a) = \alpha_i a_i$ for some $\alpha_i \in F$ and $i = 1,2,3$.
Then $a_3 = \beta_1 a_1 + \beta_2 a_2$ for some $\beta_1, \beta_2 \in F$. By expressing $f(a_3)$
in two ways, we obtain $\beta_1 \alpha_1 = \beta_1 \alpha_3$ and $\beta_2 \alpha_2 = \beta_2 \alpha_3 $.
Therefore, either $\alpha_1 = \alpha_3$ or $\alpha_2 = \alpha_3$.
\end{proof}

\medskip

\underline{\bf{2. The commutative maximal subalgebras of $M_2(\mathbb{F}_q)$.}}

\noindent Let $\mathcal{A}$ be a maximal subalgebra of $M_2({\mathbb{F}_q})$ having no nontrivial
invariant subspace in $\mathbb{F}_q^2$. Then $\mathbb{F}_q^2$ is a simple faithful
$\mathcal{A}$-module; in addition, $I_2 \in \mathcal{A}$ by Part \ref{2it} of Theorem \ref{carefree}. Therefore, by Wedderburn's Theory, $\mathcal{A}$ is isomorphic as an $\mathbb{F}_q$-algebra to a finite direct sum of matrix algebras each of which has entries in some finite extension of $\mathbb{F}_q$.
Then by the $\mathbb{F}_q$-dimension counting, $\mathcal{A}$ may be isomorphic
only to one of the following $\mathbb{F}_q$-algebras:
\begin{equation}\label{list}
\mathbb{F}_q\,,~\,\mathbb{F}_q^2\,,~\,\mathbb{F}_q^3\,,~\,\mathbb{F}_q \oplus \mathbb{F}_{q^2}\,,~\,\mathbb{F}_{q^3}\,,
~\,\mathbb{F}_{q^2} \, \text{.}
\end{equation}
The first $5$ entries in (\ref{list}) are ruled out because a simple module over a commutative ring is isomorphic to the quotient by a maximal ideal, and because of the first of Facts \ref{naive}; the resulting module has to be faithful as well. It is clear that $\mathbb{F}_q^2$
may be regarded as a simple faithful $\mathbb{F}_{q^2}$-module.
{\em It follows that up to isomorphism there is a unique maximal subalgebra $\mathcal{A}\cong \mathbb{F}_{q^2}$ of $M_2(\mathbb{F}_q)$ having no nontrivial invariant subspaces.}

The number
of such subalgebras of $M_2(\mathbb{F}_q)$ equals the number $a$ of matrices whose characteristic polynomial
is irreducible, divided by
\[b =\#(\mathbb{F}_{q^2} - \mathbb{F}_q) = q^2-q.
\] It is well known that there are exactly
$c = \left(q^2-q \right)/2$ monic irreducible quadratic polynomials
over $\mathbb{F}_q$, and there are exactly $d = q^2-q$ matrices in $M_2(\mathbb{F}_q)$
that have a given quadratic irreducible polynomial as their characteristic polynomial;
the latter is a particular case of the result of Reiner \cite{Reiner}.
{\em Hence there are exactly}
\[ \frac{a}{b}= \frac{cd }{b} =  \frac{q^2-q}{2}
\]
{\em maximal commutative subalgebras of
$M_2\left(\mathbb{F}_q\right)$}.

\bigskip

The matrices $ A_1, \ldots, A_m  \in M_2 \left( \mathbb{F}_q  \right)$ do not generate
the algebra $M_2 \left( \mathbb{F}_q  \right)$ if and only if they generate a smaller subalgebra $\mathcal{A}$, or in other words, if $\left( A_1, \ldots, A_m\right) \in \mathcal{A}^m$. Therefore, we may compute $\# G_{m,2}\left( \mathbb{F}_q \right)$ by the formula
\[
G_{m,2}\left( \mathbb{F}_q \right) =
\]
\begin{equation}\label{HomageHall}
M_2(\mathbb{F}_q)^m - \bigcup \left\{ \mathcal{A}^m \mid
~\mathcal{A}~\text{is a maximal subalgebra of}~M_2(\mathbb{F}_q) \right\}\text{.}
\end{equation}

The only $F_q$-subalgebra of $M_2(\mathbb{F}_q)$ isomorphic to $\mathbb{F}_q$
is the subalgebra of scalar matrices, which we denote by $\mathcal{D}$.
Then from the above description we see that the intersection of any two different maximal commutative subalgebras of $M_2(\mathbb{F}_q)$ is
$\mathcal{D}$. Also, since maximal subalgebras cannot be contained in each other,
it follows that the intersection of a maximal commutative subalgebra with a noncommutative one is
$\mathcal{D}$ as well.

It remains to apply the inclusion-exclusion formula to compute $\gen_{m,2}(q)$,
and this is done in
\begin{theorem}\label{th4}
Let $m \geq 2$ an integer, and $q$ a prime power, then
\begin{equation}\label{eq3}
\gen_{m,2}(q) = \frac{q^{4m-1}+q^{2m} -q^{3m}-q^{3m-1}}{q^2-1}.
\end{equation}
The function $\gen_{m,2}(q)$ is strictly increasing in $m$ and $q$. Furthermore,
we have the following equation for the generating power series:
\begin{equation}\label{genfun}
\sum_{m=2}^{\infty} \gen_{m,2}(q) z^{m-2} =
{\frac {\gen_{2,2}(q)}{ \left(1- z{q}^{2} \right) \left(1- z{q}^{3} \right)  \left( 1 -z{q}^{4} \right) }}.
\end{equation}
where $\gen_{2,2}(q) = q^4(q-1)$ according to (\ref{eq3}).
\end{theorem}

\begin{proof}
Formula (\ref{eq3}) is an arithmetic simplification of the formula
\begin{equation}\label{eq3'}
\gen_{m,2}(q) = \frac{q^{4m}-q^{m} -(q+1)\left(q^{3m}-q^m \right)+q\left(q^{2m}-q^m \right) }
{q \left(q^2-1 \right)}.
\end{equation}
Let $\mathcal{D}$ be the subalgebra of scalar matrices of
$M_2\left( \mathbb{F}_q  \right)$. According to (\ref{eq1}) and (\ref{HomageHall}), to compute $s = \#G_{m,2}(\mathbb{F}_q)$, we firstly need to remove from
$M_2(\mathbb{F}_q)^m - \mathcal{D}^m$ all the subsets of the form
$\mathcal{A}^m - \mathcal{D}^m$ for any maximal subalgebra $\mathcal{A}$ of $M_2(\mathbb{F}_q)$. This gives us $\# G_{m,2} \left( \mathbb{F}_q \right)$ which is the numerator of (\ref{eq3'}). Then we should divide $s$ by
$\# PGL_2 \left(\mathbb{F}_q \right) = q (q^2-1)$ thus obtaining (\ref{eq3'}).

According to the inclusion-exclusion formula applied to (\ref{HomageHall}), the numerator of (\ref{eq3'}) is obtained by the following two consecutive steps.

1. We subtract
from $\#\left(M_2(\mathbb{F}_q)^m - \mathcal{D}^m \right)$ the sum of the number of elements
in each of the sets of the form $\mathcal{A}^m - \mathcal{D}^m $ for any maximal subalgebra $\mathcal{A}$ of $M_2(\mathbb{F}_q)$. By the computations preceding this theorem, the resulting number is
\begin{equation}\label{eq4}
q^{4m}-q^m - (q+1)\left(q^{3m}-q^m \right) - \frac{1}{2}(q^2-q)\left(q^{2m}-q^m \right).
\end{equation}

2. Then we need to add to (\ref{eq4}) the total number of elements in each of the sets
of the form $\mathcal{A}_{u,v}^m - \mathcal{D}^m$ where
$u \neq v \in \mathbb{P}_1(\mathbb{F}_q)$. As we have seen above, this number is
$\binom{q+1}{2} \left(q^{2m}-q^m\right)$.

The application of the inclusion-exclusion formula stops at this stage,
because any intersection of $3$ pairwise different sets of the form $\mathcal{A}^m - \mathcal{D}^m$,
where $\mathcal{A}$ is a maximal subalgebra of $M_2\left( \mathbb{F}_q  \right)^2$,
is empty.

To obtain (\ref{eq3}), it remains to divide the numerator of (\ref{eq3'}) by
$\# PGL_2 \left( \mathbb{F}_q \right)$.

To see that the function $\gen_{m,2}(q)$ is strictly increasing in both variables, we write $\gen_{m,2}(q)$ as a product $f(m,q) \, g(m,q)$ of two strictly increasing functions $f(m,q) = \left(q^2 -1 \right)^{-1}q^{2m}$ and
$g(m,q) = {q}^{2\,m-1}+1-{q}^{m}-{q}^{m-1}$.

(\ref{genfun}) may be obtained from (\ref{eq3}) by using a computer algebra system.
\end{proof}

For a fixed integer $m \geq 2$, all the functions (\ref{eq3}) are polynomials in $q$ with integer coefficients, the first $3$ of which are
$\gen_{2,2}(q) = q^4(q-1)$,
$\gen_{3,2}(q) = \left( q-1 \right)  \left( {q}^{2}+q+1 \right) {q}^{6}$,
$\gen_{4,2}(q) = \left( q-1 \right)  \left( {q}^{2}+q+1 \right)
\left( {q}^{2}+1 \right) {q}^{8}$.

We note that $\gen_{2,2}(2) = 16$ according to (\ref{eq3}). The second author saw this result obtained
by a MAGMA computation by Nigel Boston \cite{Boston}, and this was a very important
piece of information that served as a guidance for this paper.

From (\ref{eq3}) it follows that
$\gen_{m,2}(q) \thicksim q^{4m-3}$ as $m,q \to \infty$. 
Therefore, it is possible that
$\gen_{m,n}(q)$ is asymptotically equivalent to (\ref{eq8}) under more general conditions than
$q \to \infty$ while $m,n$ are fixed.

\subsection{The smallest number of generators for finite direct sums of the ring $M_2(\mathbb{Z})$}\label{subs15}

Theorem \ref{th4} above yields the following estimate:
\begin{equation}\label{wow}
\gen_{m,2}(\mathbb{Z}) \leq \min_{q } \gen_{m,2}\left(q \right) =
\gen_{m,2}\left(2 \right) =
\frac{16^m - 3 \cdot 8^{m} + 2 \cdot 4^m}{6}
\end{equation}
Our goal is to show that (\ref{wow}) is an equality. This conclusion is stated
in Theorem \ref{th5} below. For $m \geq 3$, this follows from Theorem \ref{Lenstra} above. The case $m = 2$ has to be analyzed separately. We explicitly construct
$2$ generators for the ring $M_2(\mathbb{Z})^{16}$ in Table \ref{tablegen} below. We apply Theorems \ref{th1}, \ref{th3}, and Burnside's Theorem \cite{Burnside} to verify that the proposed $2$ elements generate
$M_2(\mathbb{Z})^{16}$. (Then (\ref{wow}) and Theorem \ref{gap} tell us, for example, that the smallest number of generators for $M_3(\mathbb{Z})^{17}$ is $3$.)

We will identify the group $PGL_2(\mathbb{F}_2)$ with $GL_2(\mathbb{F}_2)$ which acts on $M_2(\mathbb{F}_2)$ by conjugation, thereby creating $6$ conjugacy classes.
The two of them, those of the zero and the identity matrices, have one element each. We list the other
$4$ classes in Table \ref{tableconj} together with the corresponding eigenvalues, the eigenvalues being valid over any field. We remark that the $6$ elements of $PGL_2(\mathbb{F}_2)$ may be regarded as having entries in any
field.

\begin{table}
\begin{center}
\begin{tabular}{|c|c|}\hline
Conjugacy classes    & Eigenvalues \\ \hline
\begin{tabular}{c}
E$_{11}$,
~ E$_{22}$, ~
$
\left(
\begin{array}{cc}
1 & 0 \\
1 & 0
\end{array}
\right)$,~
$
\left(
\begin{array}{cc}
0 & 0 \\
1 & 1
\end{array}
\right)$,
\\
$\left(
\begin{array}{cc}
1 & 1 \\
0 & 0
\end{array}
\right)$,~
$\left(
\begin{array}{cc}
0 & 1 \\
0 & 1
\end{array}
\right)$
\end{tabular}
&
0, 1 \\ \hline
E$_{12}$, ~ E$_{21}$, ~
$\left(
\begin{array}{cc}
1 & 1 \\
1 & 1
\end{array}
\right)$
&
\begin{tabular}{c}
$0,0$ for the first two matrices
\\
$0,2$ for the third matrix
\end{tabular}
\\ \hline
$\left(
\begin{array}{cc}
0 & 1 \\
1 & 1
\end{array}
\right)$,~
$\left(
\begin{array}{cc}
1 & 1 \\
1 & 0
\end{array}
\right)$
&
$\text{the roots of }\lambda^2-\lambda-1 = 0$
\\ \hline
$\left(
\begin{array}{cc}
1 & 0 \\
1 & 1
\end{array}
\right)$,~
$\left(
\begin{array}{cc}
1 & 1 \\
0 & 1
\end{array}
\right)$,~
$\left(
\begin{array}{cc}
0 & 1 \\
1 & 0
\end{array}
\right)$
&
\begin{tabular}{c}
$1,1$ for the first two matrices
\\
$\pm 1$ for the third matrix
\end{tabular}
\\ \hline
\end{tabular}
\caption{Nontrivial conjugacy classes of $PGL_2(\mathbb{F}_2)$, and
the corresponding eigenvalues over any field with prime number of elements.}\label{tableconj}
\end{center}
\end{table}

We see that the eigenvalues of each of the $4$ conjugacy classes
are different modulo any prime. Therefore if we regard the above matrices as having
entries in an arbitrary prime field $\mathbb{F}_p$, then no two matrices taken from any two different
rows of Table \ref{tableconj} are conjugate modulo $p$.

The two generators for $M_2(\mathbb{Z})^{16}$ are displayed in Table \ref{tablegen}. Each entry in the
table displays the corresponding component of the two generators. According to Theorems \ref{th1},
\ref{th3}, and Burnside's Theorem \cite{Burnside}, we need to verify that

\textbullet~No pair of matrices in Table \ref{tablegen} has an eigenvector over any field
$\mathbb{F}_{p^2}$ for any prime $p$, and this is straightforward. Alternatively,
one can do an easy computer verification by Theorem \ref{modernch}.

\textbullet~No two pairs of matrices in Table \ref{tablegen} are conjugate to each
other modulo some prime. According to Table \ref{tableconj}, we only need to check this for the
two pairs marked by $\blacktriangle$ and $\blacktriangledown$ in Table \ref{tablegen}, because for any other two pairs this is automatic
by looking at the corresponding eigenvalues in Table \ref{tableconj}. This verification is straightforward and is omitted. Therefore, we have proved that the ring
$M_2\left( \mathbb{Z} \right)^{16}$ has $2$ generators, which combined with Theorem \ref{Lenstra}
gives the following
\begin{table}
\begin{center}
\begin{tabular}{|c|c|}\hline
E$_{11}$,\,
$\left(
\begin{array}{cc}
0 & 1 \\
1 & 0
\end{array}
\right)$
&
E$_{11}$, \,
$\left(
\begin{array}{cc}
1 & 1 \\
1 & 1
\end{array}
\right)$
\\ \hline
$\blacktriangle$~
E$_{11}$,\,
$\left(
\begin{array}{cc}
1 & 1 \\
1 & 0
\end{array}
\right)$
&
$\blacktriangle$~
$\text{E}_{11}, \,
\left(
\begin{array}{cc}
0 & 1 \\
1 & 1
\end{array}
\right)$
\\ \hline
E$_{12}$,~
E$_{21}$ & E$_{12}$,~
$\left(
\begin{array}{cc}
0 & 1 \\
1 & 1
\end{array}
\right)$
\\ \hline
E$_{12}$,~
$\left(
\begin{array}{cc}
0 & 1 \\
1 & 0
\end{array}
\right)$
&
E$_{12}$,~
$\left(
\begin{array}{cc}
0 & 0 \\
1 & 1
\end{array}
\right)$
\\ \hline
$\left(
\begin{array}{cc}
0 & 1 \\
1 & 1
\end{array}
\right)$,~
$\left(
\begin{array}{cc}
0 & 1 \\
1 & 0
\end{array}
\right)$
&
$\left(
\begin{array}{cc}
0 & 1 \\
1 & 1
\end{array}
\right)$,~
E$_{12}$
\\ \hline
$\blacktriangledown$~
$\left(
\begin{array}{cc}
0 & 1 \\
1 & 1
\end{array}
\right)$~,
E$_{11}$
&
$\blacktriangledown$~
$\left(
\begin{array}{cc}
0 & 1 \\
1 & 1
\end{array}
\right)$,~
E$_{22}$
\\ \hline
$\left(
\begin{array}{cc}
0 & 1 \\
1 & 0
\end{array}
\right)$,~
E$_{11}$
&
$\left(
\begin{array}{cc}
0 & 1 \\
1 & 0
\end{array}
\right)$,~
E$_{12}$
\\ \hline
$\left(
\begin{array}{cc}
0 & 1 \\
1 & 0
\end{array}
\right)$,~
$\left(
\begin{array}{cc}
0 & 1 \\
1 & 1
\end{array}
\right)$
&
$\left(
\begin{array}{cc}
0 & 1 \\
1 & 0
\end{array}
\right)$,~
$\left(
\begin{array}{cc}
1 & 1 \\
0 & 1
\end{array}
\right)$
\\ \hline
\end{tabular}
\caption{Generators for the ring $M_2(\mathbb{Z})^{16}$.}\label{tablegen}
\end{center}
\end{table}

\begin{theorem}\label{th5}
The ring $M_2(\mathbb{Z})^{16}$ has $2$ generators, while the smallest
number of generators of the ring $M_2(\mathbb{Z})^{17}$ is $3$. More generally,
\begin{equation}\label{eq152}
\gen_{m,2}(\mathbb{Z}) = \frac{16^m - 3 \cdot 8^{m} + 2 \cdot 4^m}{6}.
\end{equation}
\end{theorem}

\begin{appendices}
\section{A finite direct sum of matrix algebras over an infinite field has $2$ generators}\label{Kemer}

Theorem 3.11 (2) of \cite{PetSid} gives a presentation with $2$ generators and finitely many relations for a finite direct sum of matrix algebras over $\mathbb{Q}$.
The result almost verbatim extends to finite direct sums of matrix algebras over any infinite field $F$. The proof, however, does not if the sum has the terms of the form $M_n(F)$ where $n$ and the
characteristic of $F$ are not relatively prime. Below we prove the result for this more general case and give an application in Corollary \ref{swimming} which was not observed in \cite{PetSid}.

We recall some basic terminology used in \cite{PetSid}.
Let $F\{ x,y \}$ be the free noncommutative associative ring
whose elements we refer to as noncommutative polynomials.
For any finite direct sum $R$ of matrix algebras over $F$ we will construct an
$F$-algebra epimorphism from $F\{ x,y \}$ onto $R$.

Our considerations will be based on the following
two matrices:
\begin{equation}\label{boring_party}
X =  E_{21}+E_{32}+ \ldots + E_{n,n-1}+E_{1n} ~ ~ \text{and}~ ~
Y = E_{11} ~\text{for} ~ ~  n \geq 2 \text{.}
\end{equation}

We will use the following noncommutative polynomials:
\[
r_{1,n}  =  r_{1,n}(x) = x^n - 1, ~ \,
r_{2,n}  =  r_{2,n}(x,y) = \sum_{i=0}^{n-1} x^{n-i}yx^i - 1, \]
\[
s_0  =  s_0(y) = y^2 - y, ~ \,
s_j   =  s_j(x,y) = y x^j y ~ \, \text {for} ~ \, j \geq 1 \text{.}
\]
The matrices $X$ and $Y$ in (\ref{boring_party}) define an $F$-algebra epimorphism
from $F\{ x,y \}$ onto $M_n(F)$ by assigning $x \mapsto X$ and $y \mapsto Y$. This gives the following presentation
of the matrix algebra $M_n(F)$ introduced in \cite{PetSid}.

\begin{equation}\label{NovyyGod}
M_n(F) \cong \langle x,y \mid r_{1,n} = r_{2,n}  = s_j = 0, ~
 1 \leq j \leq n-1  \rangle.
\end{equation}

The ring $F \{ x,y \}$ has an infinite family of ideals
$\{ \mathcal{I}_n(a) \}_{a \in F}$ defined by
\[
\mathcal{I}_n (a) =
\left(r_{1,n}(x, ax + y),~ r_{2,n}(x, ax + y),~ s_j(x, ax + y),~ 1 \leq j \leq n-1  \right).
\]
Let $\mathcal{I}_n = \mathcal{I}_n (0)$. For each such an ideal, the $F$-algebras
$ M_n(F)$ and $ F \{ x,y \} / \mathcal{I}_n(a)$
are isomorphic.

\begin{theorem}\label{BraveSaid}
Let $S_1, \ldots, S_k$ be finite subsets of a field $F$, and let
$m_1, \ldots, m_k \geq 2$ be pairwise different integers, then
we have the following isomorphism of $F$-algebras:
\begin{equation}\label{Ottawa}
\frac{F \{x,y \}}{\bigcap_{i=1}^k\bigcap_{s_i \in S_i} \mathcal{I}_{m_i}(s_i)}
 \cong \bigoplus_{i=1}^k M_{m_i}(F)^{\# S_i}.
\end{equation}
\end{theorem}

\begin{proof}
The result will follow from the Chinese Remainder Theorem if we show that the sum of any two ideals in the denominator of the left hand side of (\ref{Ottawa}) is $F \{x,y \}$. Because $M_n(F)$ is a simple algebra and because of (\ref{NovyyGod}), it suffices to show that all these ideals are different, i.e.~$\mathcal{I}_m(a) = \mathcal{I}_{n}(b)$ if and only if $m = n$ and $a = b$.

\textbf{1}. Let us consider the case $m = n$. Suppose that our claim is false, so that
$
\mathcal{I}_n(a) = \mathcal{I}_n(b)
$
for some $a \neq b$.
For each $a \in F$ we introduce the $F$-algebra automorphism $\varphi_a$
of $F \{ x,y\}$ given by $\varphi_a (x) = x$ and
$\varphi_a (y) = ax + y$. We see that $\varphi_a \circ \varphi_b = \varphi_{a+b}$
and different $a$ give different $\varphi_a$.
Then
\[
\mathcal{I}_n = \varphi_{-b} \left( \mathcal{I}_n(b) \right) = \mathcal{I}_n(a-b)=
\mathcal{I}_n(c),~\text{where}~ 0 \neq c = a-b \in F.
\]
All the computations in this this first part of the proof will be done modulo the ideal $\mathcal{I}_n =\mathcal{I}_n(c)$.
\[
0 = (cx + y)x(cx + y) = c^2x^3 + c x^2y + cyx^2 + yxy = c^2x^3 + c x^2y + cyx^2.
\]
Therefore
\begin{equation}\label{smile}
0 = cx^3 +  x^2y + yx^2.
\end{equation}
If $n = 2$ then (\ref{smile}) becomes $0 = cx + 2y$. If $\text{char}\,F = 2$, then
$cx = 0$, which is impossible because both $c$ and $x$ are invertible.
If $\text{char}\,F \neq 2$, then multiplying on the left and on the right by $y$ and using $y^2 =y$ and $yxy = 0$, we obtain $y = 0$. But then $0 = r_{2,n}(x,y) = -1$,
a contradiction.

If $n \geq 3$, then multiplying (\ref{smile}) by $y$ on the right and then by
$x^{-2}$ on the left gives $0 = cx +y$, and as in the previous paragraph, we see that
$y = 0$ yielding $0 = r_{2,n}(x,y) = -1$,
a contradiction.

\bigskip

\textbf{2}. Let $m < n$ and suppose that
$
\mathcal{I}_n(a) = \mathcal{I}_n(b)
$
for some $a, b \in F$.
Then as in the previous case, we may assume that $b = 0$.
All the computations in this this second part of the proof will be done modulo the ideal $\mathcal{I}_m(a) =\mathcal{I}_n$. Then
\[
0 = y x^m y = y^2 = y \text{, so that }
0 = r_{2,n}(x,y) = -1 \text{,}
\]
a contradiction.
\end{proof}

\begin{corollary}\label{swimming}
Let $G$ be a finite group and $F$ be an algebraically closed field whose
characteristic does not divide the order of $G$. Then the group algebra
$FG$ has $2$ generators.
\end{corollary}

\begin{proof}
It is well known that the group algebra $FG$ is isomorphic to a finite direct sum
of matrix algebras over $F$, and it remains to apply Theorem \ref{BraveSaid}.
\end{proof}

A similar statement is no longer true in general for infinite groups because, for example, the group algebra $F<x,y>$ of a free Abelian group $<x,y>$ of rank $2$ does not admit two generators -- otherwise in the polynomial ring $F[u,v]$ there would exist a nonconstant invertible polynomial.


\end{appendices}

\end{document}